\documentclass[12pt]{article}
\usepackage[final]{epsfig}
\usepackage{graphics}
\usepackage{amsmath}
\usepackage{amsfonts}
\usepackage{latexsym}
\usepackage{amssymb}
\usepackage{graphicx}
\usepackage{url}
\usepackage{epstopdf}

\newtheorem{lemma}{Lemma}[section]
\newtheorem{proposition}[lemma]{Proposition}
\newtheorem{remark}[lemma]{Remark}

\newtheorem{theorem}{Theorem}

\newtheorem{corollary}[theorem]{Corollary}

\begin{document}

\newcommand{\eps}{{\varepsilon}}
\newcommand{\g}{{\gamma}}
\newcommand{\dd}{{\delta}}
\newcommand{\GG}{{\Gamma}}
\newcommand{\dD}{{\Delta}}
\newcommand{\proofend}{$\Box$\bigskip}
\newcommand{\C}{{\mathbb C}}
\newcommand{\Q}{{\mathbb Q}}
\newcommand{\R}{{\mathbb R}}
\newcommand{\Z}{{\mathbb Z}}
\newcommand{\RP}{{\mathbb {RP}}}
\newcommand{\CP}{{\mathbb {CP}}}
\newcommand{\pP}{{\mathbb {P}}}
\newcommand{\Tr}{\rm Tr}
\newcommand{\T}{\stackrel{c}{\sim}}
\newcommand{\Tt}{\stackrel{c_1}{\sim}}
\newcommand{\Ttt}{\stackrel{c_2}{\sim}}

\def\proof{\paragraph{Proof.}}

\newcommand{\SL}{\operatorname{SL}}
\newcommand{\sL}{\operatorname{sl}}
\newcommand{\PGL}{\operatorname{PSL}}
\newcommand{\Id}{\operatorname{Id}}

\newcommand{\DD}{\mathcal{D}}
\newcommand{\LL}{\mathcal{L}}
\newcommand{\CC}{\mathcal{C}}
\newcommand{\TT}{\mathcal{T}}
\newcommand{\aA}{\mathcal{A}}
\newcommand{\PP}{\mathcal{P}}

\title{On centro-affine curves and B\"acklund transformations of the KdV equation}

\author{Serge Tabachnikov\footnote{
Department of Mathematics,
Penn State University,
University Park, PA 16802;
tabachni@math.psu.edu}
}

\date{}
\maketitle

\section{A family of transformations on the space of curves}

This note stems from \cite{AFIT} where we study the integrable dynamics of a 1-parameter family of correspondences on ideal polygons in the hyperbolic plane and hyperbolic space:  two $n$-gons $P=(p_1,p_2,\ldots)$ and $Q=(q_1,q_2,\ldots)$ in $\RP^1$ or in $\CP^1$  are in correspondence  $P\T Q$  if $[p_i,p_{i+1},q_i,q_{i+1}]=c$ for all $i$; the constant $c$ is a parameter.

In the limit $n\to\infty$, a polygon becomes a  parameterized curve.
The ground field can be either $\R$ of $\C$; to fix ideas, choose $\R$. Let us
use the following definition of cross-ratio to define our correspondence (other five definitions result in the change of the constant $c$):
\begin{equation} \label{cross}
[p_i,p_{i+1},q_i,q_{i+1}]=\frac{(q_{i+1}-q_i)(p_{i+1}-p_i)}{(q_i-p_i)(q_{i+1}-p_{i+1})}=c.   
\end{equation}

We replace polygons by non-degenerate closed curves $\g: \R \to \RP^1$ with $\g'(t) >0$; to be concrete, let the period be $\pi$: $\g(t+\pi)=\g(t)$. 
Also let us assume that the rotation number of the curve $\g$ is 1, that is, $\g: \R/\pi\Z \to \RP^1$ is a diffeomorphism.
Denote the space of such curves by $\widetilde\CC$ and let $\CC=\widetilde\CC/\PGL(2,\R)$ be the moduli space. 

Then a continuous version of (\ref{cross}) is
\begin{equation} \label{Ric1}
\frac{\g'(t)\dd'(t)}{(\dd(t)-\g(t))^2}=c.
\end{equation}

Write $\g \T \dd$ to denote this relation on $\widetilde\CC$. Since cross-ratio is M\"obius invariant, we also have a relation on $\CC$ which we continue to denote by $\T$. Note that $\T$ is a symmetric relation.

\begin{lemma} 
 For a generic curve $\g$, the relation $\g \T \dd$ is a (partially defined) 2-2 map $T_c: \g \mapsto \dd$. 
\end{lemma}

\proof
Given $\g(t)$, equation (\ref{Ric1}) is a Riccati equation on $\dd(t)$, its monodromy is a M\"obius transformation (see, e.g., \cite{In}) which has either two or no fixed points, unless it is the identity. Over $\C$, there are always two fixed points (possibly, coinciding), and over $\R$, we need to assume that they exist. Then $\T$ defines a 2-2 map.
\proofend

Thus, given $\g$, there are two choices of $\dd=T_c (\g)$. Once a choice is made, one similarly has two choices for $T_c(\dd)$, but one of them is $\g$, so we choose the other one, and so on. Hence the choice of $\dd$ determines the map $T_c$; the other choice gives the inverse map $T_c^{-1}$.

Following the standard procedure (see, e.g., \cite{OT}),  lift a curve $\g(t)$ from $\RP^1$ to $\R^2$, normalizing the lift $\GG(t)$ so that $[\GG,\GG']=1$ (here and elsewhere $[\, , ]$ denotes the determinant made by two vectors). 

Explicitly, $\GG=((\g')^{-1/2}, (\g')^{-1/2} \g)$. Note the  square root: the curve $-\GG$ will do as well, the lift is defined up to the sign, and the action of $\PGL(2,\R)$ is replaced by that of $\SL(2,\R)$.
 We obtain centro-affine realizations of the spaces $\widetilde\CC$ and $\CC$. 

The curve $\GG$ satisfies a Hill equation
\begin{equation} \label{Hill}
\GG''(t)=p(t)\GG(t)
\end{equation}
with a $\pi$-periodic potential $p(t)$, and $\GG(t+\pi)=-\GG(t)$ (the curve makes exactly half-rotation on $[0,\pi]$). In geometric terms, the potential $p$ is the (negative) centro-affine curvature of the curve $\GG$.

In these terms, equation (\ref{cross})  becomes
$$
\frac{[\dD(t),\GG(t)][\dD(t+\eps),\GG(t+\eps)]}{[\dD(t),\dD(t+\eps)][\GG(t),\GG(t+\eps)]} = {\rm const},
$$
and, in the limit $\eps\to 0$, we obtain an analog of equation (\ref{Ric1}):
$$
[\GG(t),\dD(t)]^2=c^2.
$$
Break the symmetry between $\GG$ and $\dD$ by taking square root:
\begin{equation} \label{liftdef}
[\GG(t),\dD(t)]=c.
\end{equation}
This defines a map on the lifted curves: $T_c(\GG)=\dD$. Note that $T_c(\dD)=-\GG$.

\begin{lemma} \label{twomap}
$T_c: \GG \mapsto \dD$ is a (partially defined) 2-2 map. 
\end{lemma}

\proof
Let us search for $\dD$ in the form $\dD = a \GG + b \GG'$, where $a$ and $b$ are $\pi$-periodic functions. Then equation $[\GG,\dD]=c$ implies that $b(t)=c$, $\dD=a\GG+c\GG'$, and then $\dD'=(a'+cp)\GG + a\GG'$. The condition $[\dD,\dD']=1$ now implies 
\begin{equation} \label{Ric}
a' = \frac{a^2-1}{c} - cp.
\end{equation}
This is a Riccati equation on function $a(t)$ with periodic coefficients. The monodromy of this equation is a M\"obius transformation, hence it has two fixed points (always, if one works over $\C$, and over $\R$ one needs to assume that it does), corresponding to two periodic solutions of (\ref{Ric}). Each solution defines a curve $\dD$ with $T_c(\GG)= \dD$.
\proofend

As before, once a choice of a fixed point of the monodromy is made, the map becomes 1-1: of the two choices available for the next curve $\dD$, one is extraneous because it takes one  back to $-\GG$. 

\section{Two pre-symplectic forms and a bi-Hamiltonian structure}

Starting with U. Pinkall \cite{Pin}, a number of recent papers were devoted to the study of the Korteweg-de Vries equation in terms of cento-affine curves \cite{CIMB,FK1,FK,TW}. Let us present the relevant results.

A tangent vector at a cento-affine curve $\GG$ is a vector field along $\GG$ that can be written as a linear combination $h\GG+f\GG'$ where $h,f$ are $\pi$-periodic functions.

\begin{lemma}
The function $f$ is arbitrary, and 
$h=-\frac{1}{2} f'.$
\end{lemma}

\proof
If $\GG+\eps v$ is a deformation of $\GG$, then $[\GG,v']=[\GG',v]$ because $[\GG,\GG']=1$. For $v=h\GG+f\GG'$, this implies that $h=-\frac{1}{2} f'.$
\proofend

Denote the tangent vectors by $U,V$ or $U_f,V_f$, in the format $-\frac{1}{2} f' \GG + f \GG'$.

The following pre-symplectic structure on space $\widetilde\CC$ was introduced in \cite{Pin}.
Let $U,V$ be tangent vector fields along $\GG$; define
$$
\omega(U,V) = \int_\GG [U,V]\ dt,
$$
that is,
$$
\omega(U_f,V_g) = \frac{1}{2}\int_0^{\pi} (fg'-f'g)dt.
$$
The kernel of $\omega$ is spanned by the field $\GG'$, that is, by the reparameterizations $t \mapsto t+$const.

Pinkall observed that the Hamiltonian vector field of the function $\int p\ dt$ is $U_p$, which induces the KdV evolution 
of the potential $p$ 
$$
\dot p = -\frac{1}{2} p''' + 3p'p
$$
(the signs differ from those of Pinkall because he used the opposite sign for the potential of Hill's equation).

The second pre-symplectic structure was introduced in \cite{FK}: for tangent vector fields $U,V$  along $\GG$, let
$$
\Omega(U,V) =  \int_\GG ([U',V'] +p[U,V])\ dt,
$$
that is,
$$
\Omega(U_f,V_g) = \int_0^{\pi} \left[\frac{1}{4}(f'g''-f''g') + p(fg'-f'g)\right]\ dt.
$$

Concerning the kernel of $\Omega$, one has 

\begin{lemma} [\cite{FK}] \label{kerOmega}
The kernel of $\Omega$ is 3-dimensional, it is generated by the Killing vector fields $A(\GG)$ with $A \in \SL(2)$.
\end{lemma}

\proof
One has
$$
\Omega(U,V) = \int [pU-U'',V]\ dt.
$$
Hence $U$ is in the kernel if and only if $U''=pU$, that is, $U(t)$ is $\SL(2)$-equivalent to $\GG(t)$. 
\proofend

Thus the form $\Omega$ descends on the moduli space $\CC$ as a symplectic form.

It is shown in \cite{CIMB,TW,FK} that the forms $\omega$ and $\Omega$ provide a bi-Hamiltonian structure on $\widetilde\CC$, corresponding to a pair of compatible Poisson brackets for the KdV equation.

Namely, let $X_0,X_1,\ldots$ and $H_1,H_2,\ldots$ be the vector fields and the Hamiltonians of the KdV hierarchy in terms of centro-affine curves:
$$
X_0=U_1=\GG', X_1=U_p=-\frac{p'}{2} \GG + p \GG',\ldots , H_1 = \int p dt, H_2 = \frac{1}{2} \int p^2 dt, \ldots
$$
Then one has
\begin{equation} \label{reccur}
\Omega (X_{j-1},\cdot) = d H_j = \omega(X_j,\cdot),\ j=1,2,\ldots
\end{equation}
see \cite{FK}.
\paragraph{The forms $\omega$ and $\Omega$ on projective curves.}
Let us calculate these forms in terms of the curves $\g: \R \to \RP^1$.

 In \cite{AFIT}, the following differential 2-form on the space of polygons $(p_1,\ldots,p_n) \subset \RP^1$ was considered
$$
\omega'=\sum_i \frac{dp_i \wedge dp_{i+1}}{(p_{i+1}-p_i)^2},
$$
and it was proved that this form was $T_c$-invariant. In the continuous limit, a polygon becomes a curve $\g(t)$.
Let $u(t),v(t)$ be two vector fields along $\g(t)$, that is, two periodic functions. Then, in the continuous limit,  we obtain the form
$$
\omega'(u,v) = \int \frac{uv'-u'v}{(\g')^2}\ dt.
$$

\begin{lemma} \label{twoomegas}
One has $\omega=\frac{1}{2} \omega'.$
\end{lemma}

\proof
Since
$$
\GG_1=(\g')^{-1/2},\ \GG_2=(\g')^{-1/2} \g,
$$
one calculates the respective vector field along $\GG$:
$$
U = \left( -\frac{1}{2} u' \GG_1^3, -\frac{1}{2} u' \GG_1^2\GG_2+u\GG_1  \right),
$$
and likewise for $V$. Then
$$
[U,V]=\frac{1}{2} \GG_1^4 (uv'-u'v),
$$
and the result follows.
\proofend

By Lemma \ref{kerOmega}, the 2-form $\Omega$  descends to the moduli space of projective curves, that is, to the space of Hill's equations. This space is 
a coadjoint orbit of the Virasoro algebra, and $\Omega$ coincides (up to a factor) with the celebrated Kirillov-Kostant-Souriau symplectic structure, see, e.g., \cite{KW,OT} for this material.  

Namely, let $\g$ be a curve in $\RP^1$, and let $u$ and $v$ be vector fields along $\g$.
The Kirillov-Kostant-Souriau symplectic form is given by the formula
$$
\Omega'(u,v)=\int \frac{u''(t)v'(t)-u'(t)v''(t)}{\g'(t)^2}\ dt,
$$
see, e.g., \cite{OT1}. 

\begin{lemma}
One has
$$
\Omega = -\frac{1}{4} \Omega'.
$$
\end{lemma} 

\proof
As in the proof of Lemma \ref{twoomegas},
$$
U = \left( -\frac{1}{2} u' \GG_1^3, -\frac{1}{2} u' \GG_1^2\GG_2+u\GG_1  \right),
$$
and then
$$
U'=\left( -\frac{1}{2} u''\GG_1^3-\frac{3}{2} u'\GG_1^2\GG_1', -\frac{1}{2} u''\GG_1^2\GG_2 - u'\GG_1\GG_2\GG_1' -\frac{1}{2} u'\GG_1^2\GG_2'+u'\GG_1+u\GG_1' \right).
$$
Similar formulas hold for $V$.

Now one computes, using the fact that $\GG_2'\GG_1-\GG_1'\GG_2=1$, 
$$
[U',V'] = -\frac{1}{4}\GG_1^4(u''v'-u'v'') - \frac{1}{2} \GG_1^3\GG_1' (u''v-uv'') - \frac{3}{2} \GG_1^2(\GG_1')^2 (u'v-uv'),
$$
and
$$
p[U,V] = -\frac{1}{2} p\GG_1^4 (u'v-uv') = -\frac{1}{2} \GG_1^3 \GG_1'' (u'v-uv').
$$
Integrating by parts,
$$
- \int \GG_1^3\GG_1' (u''v-uv'')\ dt = \int (\GG_1^3\GG_1'' + 3 \GG_1^2 (\GG_1')^2 ) (u'v-uv')\ dt,
$$
and collecting terms, 
$$
\Omega(U,V) = \int ([U',V']+p[U,V])\ dt = -\frac{1}{4} \int \GG_1^4(u''v'-u'v'')\ dt = -\frac{1}{4} \Omega'(u,v),
$$
as claimed.
\proofend

\section{$T_c$-invariance of the bi-Hamiltonian structure and complete integrability of the transformations  $T_c$}

Let $T_c(\GG)=\dD$ with $\dD''=q \dD$; one can write $\dD(t)=a(t) \GG(t) + c \GG'(t)$, where $a(t)$ is a periodic function.

\begin{lemma} \label{rels}
One has: 
$$
\GG=a\dD-c\dD',\ p+q=\frac{2}{c^2} (a^2-1),\ q-p=\frac{2}{c} a'.
$$
\end{lemma}

\proof
Since $[\dD,-\GG]=c$, we can write $-\GG=b\dD+c\dD'$ where $b(t)$ is a periodic function. Substitute $\dD=a\GG+c\GG'$ in this equation to find that $b=-a$. We also have an analog of (\ref{Ric}) for function $b(t)$:
$cb'=b^2-1-c^2q$. This implies the relations between $p$ and $q$ stated in the lemma.
\proofend

Let $T_c(\GG)=\dD$, and let $U_f, V_g$
be two tangent vectors, at $\GG$ and $\dD$, respectively,  related by the differential of $T_c$.

\begin{lemma} \label{differ}
One has
\begin{equation} \label{refvect}
\frac{c}{2} (f'+g') = a(g-f),
\end{equation}
where the function $a(t)$ is as above.
\end{lemma} 

\proof
One has $[U,\dD]+[\GG,V]=0$, or
$$
\frac{c}{2} (f'+g') = f[\GG',\dD]+g[\GG,\dD'] = a(g-f),
$$
where the last equality makes use of $\dD = a \GG + c \GG'$ and of $[\GG',\dD]+[\GG,\dD']=0$.
\proofend

The following theorem is our main observation.

\begin{theorem} \label{invomega}
The forms $\omega$ and $\Omega$ are invariant under the maps $T_c$: 
$$
T_c^*(\omega)=\omega,\ T_c^*(\Omega)=\Omega.
$$
\end{theorem}

\proof
Let $T_c(\GG)=\dD$, and let $U_{f_i}, V_{g_i}, i=1,2$, be two pairs of tangent vectors, at $\GG$ and $\dD$, respectively, related by the differential of $T_c$.
One has
\begin{equation*}
\begin{split}
\int [(g_1' g_2-g_1g_2') - (f_1' f_2-f_1f_2')] dt 
= 
\int [(g_1' g_2-g_1g_2') - (f_1' f_2-f_1f_2')& \\
-  (g_1'f_2+g_1f_2') + (g_2'f_1+g_2f_1')]dt&\\
=       
\int [(f_1'+g_1')(g_2-f_2)
- (f_2'+g_2')(g_1-f_1)]dt =0&,
\end{split}
\end{equation*}
where the first equality follows from the fact that $g_1'f_2+g_1f_2'=(g_1f_2)'$ and $g_2'f_1+g_2f_1'=(g_2f_1)'$, which integrates to zero, and the last equality follows from (\ref{refvect}). Thus  $T_c^*(\omega)=\omega$.

To prove that $T_c^*(\Omega)=\Omega$, we argue similarly, although the computation is more involved.
 
Differentiate (\ref{refvect}) to obtain
\begin{equation} \label{diffref}
\frac{c}{2} (f''+g'')=a'(g-f)+a(g'-f').
\end{equation}

We want to show that the integral
\begin{equation} \label{intzero}
\int \left(\frac{1}{4} (f_1'f_2''-f_1''f_2') + p(f_1f_2'-f_1'f_2) -  \frac{1}{4} (g_1'g_2''-g_1''g_2') -q(g_1g_2'-g_1'g_2)\right) dt 
\end{equation}
vanishes. One has
$$
f_1'f_2''-f_1''f_2' -g_1'g_2''+ g_1''g_2' = (f_1''+g_1'')(g_2'-f_2') - (f_2''+g_2'')(g_1'-f_1') +(f_2'g_1'-f_1'g_2')',
$$
hence
\begin{equation*}
\begin{split}
&\frac{1}{4} \int (f_1'f_2''-f_1''f_2' -g_1'g_2''+ g_1''g_2') dt \\
= &\frac{1}{2c} \int \{[a'(g_1-f_1)+a(g_1'-f_1')](g_2'-f_2') - [a'(g_2-f_2)+a(g_2'-f_2')](g_1'-f_1')\} dt\\
= &\int \frac{a'}{2c} [(g_1-f_1)(g_2'-f_2') - (g_2-f_2)(g_1'-f_1')] dt,
\end{split}
\end{equation*}
where the first equality follows from (\ref{diffref}). 

Next we evaluate the remaining part of the integral (\ref{intzero}), using Lemma \ref{rels}:
\begin{equation*}
\begin{split}
&\int [p(f_1f_2'-f_1'f_2)-q(g_1g_2'-g_1'g_2)] dt \\
= &\int  \frac{a^2-1}{c^2}(f_1f_2'-f_1'f_2-g_1g_2'+g_1'g_2) dt  - \int  \frac{a'}{c} (f_1f_2-f_1'f_2+g_1g_2'-g_1'g_2) dt.
\end{split}
\end{equation*}
Collecting the integrals together, we obtain
\begin{equation*}
\begin{split}
&\int  \frac{a'}{2c} [(f_1'+g_1')(f_2+g_2) - (f_2'+g_2')(f_1+g_1)] dt \\
+ &\int  \frac{a^2-1}{c^2}(f_1f_2'-f_1'f_2-g_1g_2'+g_1'g_2) dt \\
= &\int  \frac{2aa'}{c^2} (f_2g_1-f_1g_2) dt + \int  \frac{a^2-1}{c^2}(f_1f_2'-f_1'f_2-g_1g_2'+g_1'g_2) dt, 
\end{split}
\end{equation*}
where the equality is due to (\ref{refvect}).

Finally, notice that $(a^2-1)'=2aa'$, and integrate by parts to obtain
\begin{equation*}
\begin{split}
\int \frac{a^2-1}{c^2} [(f_1f_2'-f_1'f_2-g_1g_2'+g_1'g_2) - (f_2'g_1 + f_2g_1'-f_1'g_2-f_1g_2')] dt \\
= \int \frac{a^2-1}{c^2} [(f_2'+g_2')(f_1-g_1)-(f_1'+g_1')(f_2-g_2)] dt =0,
\end{split}
\end{equation*}
since the last integrand vanishes due to (\ref{refvect}).
\proofend

\begin{corollary}
The maps $T_c$ commute with the KdV flows and preserve the KdV integrals. 
\end{corollary}

\proof
One argues inductively using formulas (\ref{reccur}):
$$
\Omega (X_{j-1},\cdot) = d H_j = \omega(X_j,\cdot).
$$
 If $T_c$ preserves $X_{j-1}$ then, since it also preserves $\Omega$, it preserves $dH_j$. If $T_c$ preserves $dH_j$ then, since it preserves $\omega$, it also preserves $X_j$. 
 
 To start the induction, we check that $\int p\ dt$ is invariant:
 $$
\int (q(t)-p(t))\ dt = \frac{2}{c} \int a'(t)\ dt =0
$$
due to Lemma \ref{rels}.

Since $dH_j$ is preserved, it could be that $T_c$ changes $H_j$ by a constant. To see that this constant is zero, let $\GG$ be the circle $(\cos t, \sin t)$. Then $\dD$ differs from $\GG$ by a parameter shift, and the values of the  functions $H_j$ on $\GG$ and $\dD$ are equal.
\proofend

Thus the transformations $T_c$ are symmetries of the Korteweg-de Vries equation.

\begin{remark}
{\rm 
The argument above is similar to the one given in \cite{Tab} which concerned with the filament equation and the bicycle transformations as its symmetries.
}
\end{remark}

\paragraph{Additional  integrals.}
Let $\GG=(\GG_1,\GG_2)$. Consider the functions
\begin{equation*} \label{upint}
I= \int \GG_1^2\ dt,\ J=\int \GG_1\GG_2\ dt,\ K= \int \GG_2^2\ dt
\end{equation*}
on the space of centro-affine curves.

\begin{proposition} \label{slHam}
The functions $I,J,K$ are the Hamiltonians of the generator of the action of $\sL(2,\R)$ on $\widetilde\CC$ with respect to the 2-form $\omega$.
The function $IK-J^2$ is $\SL(2,\R)$-invariant.
\end{proposition}

\proof
The generators of $\sL(2,\R)$ are the fields 
$$
(\GG_2,0),\ (\GG_1, -\GG_2),\ (0,\GG_1).
$$
Let us consider the first one; the other ones are dealt with similarly. 

We claim that
$
(\GG_2,0) = - V_{\GG_2^2}.
$
Indeed, 
$$
V_{\GG_2^2} = -\GG_2\GG_2' \GG + \GG_2^2 \GG'.
$$
The first component of this vector is $-\GG_2 (\GG_2'\GG_1 - \GG_1'\GG_2) = -\GG_2$, and the second component is $-\GG_2'\GG_2^2 + \GG_2'\GG_2^2=0$.

Let $U_f$ be a test vector field. Then 
$$
dK(U_f) = \int \GG_2\left( \GG_2'f - \frac{1}{2} \GG_2 f'\right)\ dt = 2 \int \GG_2\GG_2'f\ dt.
$$
On the other hand, 
$$
\omega(U_f,U_{\GG_2^2}) = \int 2\GG_2\GG_2'f\ dt,
$$
as needed.

As to $\sL(2,\R)$ invariance of $IK-J^2$, let us again check invariance under the field $(\GG_2,0)$ (the rest is similar). Calculating mod $\eps^2$, one has
\begin{equation*}
\begin{split}
\left(\int (\GG_1+\eps\GG_2)^2\ dt\right) \left(\int \GG_2^2\ dt\right) - \left( \int (\GG_1+\eps \GG_2)\GG_2\ dt \right)^2
= IK-J^2& \\
+ 2 \eps \left[ \left(\int \GG_1\GG_2\ dt\right) \left(\int \GG_2^2\ dt\right) - \left(\int \GG_1\GG_2\ dt\right) \left(\int \GG_2^2\ dt\right)\right] 
= IK-J^2&,
\end{split}
\end{equation*}
as needed
\proofend

Next we show that $I,J,K$ are integrals of the transformations $T_c$.

\begin{theorem}
Let $T_c(\GG)= \dD$, then 
$$
I(\GG)=I(\dD), J(\GG)=J(\dD), K(\GG)=K(\dD).
$$
\end{theorem} 

\proof Consider the case of $I$; the other two cases are similar.

We have $\dD=a\GG+c\GG'$, and we want to show that $\int \dD_1^2 = \int \GG_1^2$. Indeed,
\begin{equation*} 
\begin{split}
&\int (\dD_1^2 - \GG_1^2)\ dt = \int [(a^2-1)\GG_1^2+2ca\GG_1\GG_1'+c^2(\GG_1')^2]\ dt\\
&= \int [(a^2-1-ca') \GG_1^2 + c^2(\GG_1')^2]\ dt = c^2 \int [p \GG_1^2 + (\GG_1')^2]\ dt \\
&= c^2 \int [\GG_1''\GG_1 + (\GG_1')^2]\ dt = 0,
\end{split}
\end{equation*}
where the second equality is integration by parts, the third is due to (\ref{Ric}),  the fourth is due to $\GG''=p\GG$, and the last one is again integration by parts.
\proofend

\section{Monodromy integrals and permutability}

Now we describe an infinite collection of $\SL(2,\R)$-invariant integrals of the maps $T_c$  that arise from the monodromy of the Riccati equations. 

Let $x$ be an affine coordinate on $\RP^1$. The Lie algebra $\sL(2,\R)$ is generated by the vector fields
$ \partial_x, x \partial_x, x^2 \partial_x.$ Introduce time-dependent vector fields, depending on $\g(t)$ or $\dd(t)$, respectively, taking values in $\sL(2,\R)$ for each $t$:
$$
\xi_\g = \left(\frac{\g^2}{\g'} -2\frac{\g}{\g'}x + \frac{1}{\g'} x^2  \right) \partial_x,\ \ \xi_\dd = \left(\frac{\dd^2}{\dd'} -2\frac{\dd}{\dd'}x + \frac{1}{\dd'} x^2  \right) \partial_x.
$$
Then equation (\ref{Ric1}) describes $\dd$ as evolving under the field $c\xi_\g$ and, equivalently, $\g$ as evolving under $c\xi_\dd$.

Fix a (spectral) parameter $\lambda$, and consider the time-$\pi$ flows of the fields $\lambda\xi_\g$ and $\lambda\xi_\dd$,  where $\g$ and $\dd$ are related by (\ref{Ric1}). Denote these  projective transformations of $\RP^1$ by $\Phi_{\lambda,\g}$ and $\Phi_{\lambda,\dd}$.

\begin{theorem} \label{famint}
For every $\lambda$, the maps $\Phi_{\lambda,\g}$ and $\Phi_{\lambda,\dd}$ are conjugate in $\PGL(2,\R)$.
\end{theorem} 

It follows that the spectral invariants of $\Phi_{\lambda,\g}$, say $\Tr^2/\det$, as functions of $\lambda$, are integrals of the maps $T_c$ for all values of $c$.

\proof
Let $\g$ and $\dd$ satisfy (\ref{Ric1}). Introduce a time-dependent matrix, also depending on parameter $\mu$:
$$
A_{\mu,\g,\dd} (t) = \frac{1}{\g(t)-\dd(t)}
\begin{bmatrix}
	\g(t)-\mu\dd(t),&\g(t)\dd(t)(\mu-1)\\
	1-\mu,&\g(t)\mu-\dd(t)
	\end{bmatrix}.
$$
We claim that if $\lambda = c(1-\mu)$, then $A_{\mu,\g,\dd} (t)$ conjugates the  vector fields $\lambda\xi_\g$ and $\lambda\xi_\dd$.

Namely, let $\eps$ be an infinitesimal parameter, and set
$$
V_\g (t,\eps) = 
\begin{bmatrix}
	1-  \frac{\eps \lambda\g(t)}{\g(t)'},& \frac{\eps \lambda\g(t)^2}{\g(t)'}\\
	-\frac{\eps\lambda}{\g(t)'},&1+\frac{\eps\lambda\g(t)}{\g(t)'}
	\end{bmatrix}
$$
This time-dependent M\"obius transformation is the time-$\eps$ flow of the vector field $\lambda\xi_\g$. 

Then one has
$$
V_\dd (t,-\eps)  A_{\mu,\g,\dd} (t+\eps)  V_\g (t,\eps) = V_\dd (t,\eps)  A_{\mu,\g,\dd} (t-\eps)  V_\g (t,- \eps) \mod \eps^2,
$$
which is verified by a direct calculation or, in the limit $\eps\to 0$, 
$$
\begin{bmatrix}
	 \frac{\dd(t)}{\dd'(t)},& \frac{\dd(t)^2}{\dd'(t)}\\
	-\frac{1}{\dd'(t)},&\frac{\dd(t)}{\dd'(t)}
	\end{bmatrix} A_{\mu,\g,\dd}(t)
- A_{\mu,\g,\dd}(t) 	\begin{bmatrix}
	 \frac{\g(t)}{\g'(t)},& \frac{\g(t)^2}{\g'(t)}\\
	-\frac{1}{\g'(t)},&\frac{\g(t)}{\g'(t)}
	\end{bmatrix}
	       = \frac{1}{\lambda} A_{\mu,\g,\dd}' (t).
$$
This equality implies that the vector fields $\lambda\xi_\g$ and $\lambda\xi_\dd$ are conjugate, and so are $\Phi_{\lambda,\g}$ and $\Phi_{\lambda,\dd}$:
\begin{equation} \label{conj}
\Phi_{\lambda,\dd} = A_{\mu,\g,\dd} (0) \Phi_{\lambda,\g} A_{\mu,\g,\dd}^{-1} (0),
\end{equation}
as needed.
\proofend

\begin{remark}
{\rm The above theorem is also a continuous analog of a result for ideal polygons in \cite{AFIT}.
}
\end{remark}

\paragraph{Bianchi permutability.} Let us show that the maps $T_c$ commute; the argument is similar to that given in \cite{AFIT} for ideal polygons.

\begin{theorem} \label{permut}
Let three closed curves satisfy $\g \Tt \g_1$ and $\g \Ttt \g_2$. Then there exists a fourth curve $\g_{12}$ such that $\g_1 \Ttt \g_{12}$ and $\g_2 \Tt\g_{12}$.
\end{theorem}

\proof
We use (\ref{conj}), writing $A$ instead of $A(0)$.

Since $\g \Tt \g_1$ and $\g \Ttt \g_2$, we have
$$
\Phi_{c_1,\g} (\g_1(0)) = \g_1(0),\ \Phi_{c_2,\g} (\g_2(0)) = \g_2(0).
$$ 
By (\ref{conj}), 
$$
\Phi_{c_1,\g_2} = A_{\mu,\g,\g_2} \Phi_{c_1,\g} A_{\mu,\g,\g_2}^{-1}, \ \Phi_{c_2,\g_1} = A_{\nu,\g,\g_1} \Phi_{c_2,\g} A_{\nu,\g,\g_1}^{-1}
$$
with 
\begin{equation} \label{relc}
c_1=c_2(1-\mu), c_2=c_1(1-\nu).
\end{equation}
It follows that
$$
\Phi_{c_1,\g_2} (A_{\mu,\g,\g_2} (\g_1(0))) = A_{\mu,\g,\g_2} (\g_1(0)),\ \Phi_{c_2,\g_1} (A_{\nu,\g,\g_1} (\g_2(0))) = A_{\nu,\g,\g_1} (\g_2(0)).
$$
Thus we need to show that 
\begin{equation} \label{match}
A_{\mu,\g,\g_2} (\g_1(0)) = A_{\nu,\g,\g_1} (\g_2(0)).
\end{equation}

This is indeed the case: (\ref{relc}) implies that $\frac{1}{\mu} + \frac{1}{\nu} =1$, and then a calculation shows that
$$
\frac{1}{\mu}
\begin{bmatrix}
	\g-\mu\g_2,&\g\g_2(\mu-1)\\
	1-\mu,&\g\mu-\g_2
	\end{bmatrix}
\begin{bmatrix}
	\g_1\\
	1
	\end{bmatrix}
	=
	\frac{1}{\nu}
\begin{bmatrix}
	\g-\nu\g_1,&\g\g_1(\nu-1)\\
	1-\nu,&\g\nu-\g_1
	\end{bmatrix}
\begin{bmatrix}
	\g_2\\
	1
	\end{bmatrix},	
$$
as needed.
\proofend

\begin{remark}
{\rm The above considerations can be extended to centro-affine {\it twisted} curves, that is, curves with monodromy, $\GG(t+\pi)=M(\GG(t))$, where the monodromy $M\in\SL(2,\R)$ is not necessarily $-\Id$.  One can define the maps $T_c$ on twisted curves: given $\GG$, consider the respective $\pi$-periodic potential of the Hill equation $p(t)$, find a $\pi$-periodic solution $a(t)$ to equation (\ref{Ric}), and define $\dD=a\GG+c\GG'$. Then the monodromy of  $\dD$ coincides with that of $\GG$. At the level of Hill's equations, this is the map $p\mapsto q$. We do not dwell on this extension here.
}
\end{remark}

\bigskip
{\bf Acknowledgements}. It is a pleasure to acknowledge the stimulating discussions with A. Calini, A. Izosimov, I. Izmestiev, B. Khesin, and V. Ovsienko. This work was supported by NSF grant DMS-1510055.

\end{document}